\documentclass[12pt]{article}
\usepackage[english,frenchb]{babel}
\usepackage[latin1]{inputenc}
\usepackage{amsfonts}                                             
\usepackage{amsmath}
\usepackage{amsthm}
\usepackage{fullpage}

\begin{document}

\newcommand{\ABS}[1]{{\left| #1 \right|}} % |1|
\newcommand{\PAR}[1]{{\left(#1\right)}} % (1)
\newcommand{\SBRA}[1]{{\left[#1\right]}} % [1]
\newcommand{\BRA}[1]{{\left\{#1\right\}}} % {1}
\newcommand{\NRM}[1]{{\Vert #1\Vert}} % ||1|
\newcommand{\ind}{\mathrm{1}\hskip -3.2pt \mathrm{I}}

% Environnements etc.
\newtheorem{ethm}{Theorem}[section]
\newtheorem{fthm}{Théorème}[section]
\newtheorem{ecor}[ethm]{Corollary}
\newtheorem{fcor}[fthm]{Corollaire}
\newtheorem{eprop}[ethm]{Proposition}
\newtheorem{fprop}[fthm]{Proposition}
\newtheorem{elem}[ethm]{Lemma}
\newtheorem{flem}[fthm]{Lemme}
\newtheorem{edefi}[ethm]{Definition}
\newtheorem{fdefi}[fthm]{Définition}
\newtheorem{erem}[ethm]{Remark}
\newtheorem{frem}[fthm]{Remarque}
\newtheorem{fex}{Exemple}
\newtheorem{eex}{Example}
\newtheorem{fhyp}[fthm]{Hypothèse}

% LETTRES GRECQUES
\newcommand{\al}{\alpha}
\newcommand{\be}{\beta}
\newcommand{\de}{\delta}
\newcommand{\De}{\Delta}
\newcommand{\ga}{\gamma}
\newcommand{\Ga}{\Gamma}
\newcommand{\ep}{\epsilon}
\newcommand{\va}{\varphi}
\newcommand{\ka}{\kappa}
\newcommand{\la}{\lambda}
\newcommand{\La}{\Lambda}
\newcommand{\te}{\theta}
\newcommand{\Te}{\Theta}
\newcommand{\om}{\omega}
\newcommand{\Om}{\Omega}
\newcommand{\si}{\sigma} 
\newcommand{\ot}{\otimes}
\newcommand{\ti}{\times}
\newcommand{\cd}{CD(\rho ,\infty)}
\newcommand{\varep}{\varepsilon}
\newcommand{\na}{\nabla}
\newcommand{\ph}{\Phi}

% Double letters
\newcommand{\dA}{\mathbb{A}}
\newcommand{\dB}{\mathbb{B}}
\newcommand{\dC}{\mathbb{C}}
\newcommand{\dD}{\mathbb{D}}
\newcommand{\dE}{\mathbb{E}}
\newcommand{\dF}{\mathbb{F}}
\newcommand{\dG}{\mathbb{G}}
\newcommand{\dH}{\mathbb{H}}
\newcommand{\dI}{\mathbb{I}}
\newcommand{\dJ}{\mathbb{J}}
\newcommand{\dK}{\mathbb{K}}
\newcommand{\dL}{\mathbb{L}}
\newcommand{\dM}{\mathbb{M}}
\newcommand{\dN}{\mathbb{N}}
\newcommand{\dO}{\mathbb{O}}
\newcommand{\dP}{\mathbb{P}}
\newcommand{\dQ}{\mathbb{Q}}
\newcommand{\dR}{\mathbb{R}}
\newcommand{\dS}{\mathbb{S}}
\newcommand{\dT}{\mathbb{T}}
\newcommand{\dU}{\mathbb{U}}
\newcommand{\dV}{\mathbb{V}}
\newcommand{\dW}{\mathbb{W}}
\newcommand{\dX}{\mathbb{X}}
\newcommand{\dY}{\mathbb{Y}}
\newcommand{\dZ}{\mathbb{Z}}

% CALLIGRAPHIE
\newcommand{\cA}{{\mathcal A }}
\newcommand{\cB}{{\mathcal B }}
\newcommand{\cC}{{\mathcal C }}
\newcommand{\cD}{{\mathcal D }}
\newcommand{\cE}{{\mathcal E }}
\newcommand{\cF}{{\mathcal F }}
\newcommand{\cG}{{\mathcal G }}
\newcommand{\cH}{{\mathcal H }}
\newcommand{\cI}{{\mathcal I }}
\newcommand{\cJ}{{\mathcal J }}
\newcommand{\cK}{{\mathcal K }}
\newcommand{\cL}{{\mathcal L }}
\newcommand{\cM}{{\mathcal M }}
\newcommand{\cN}{{\mathcal N }}
\newcommand{\cO}{{\mathcal O }}
\newcommand{\cP}{{\mathcal P }}
\newcommand{\cZ}{{\mathcal Z }}

%%% Théorème sans numérotation
\newenvironment{ithm}
 {\noindent \\ \noindent{\emph{\textbf{Théorème.}}} \begin{it}}
 {\end{it}}            
\newenvironment{iprop}
 {\noindent {\emph{\textbf{Proposition.}}} \begin{it}}
 {\end{it}\\}     

%%% Preuves
\newcommand{\proofend}{\hfill $\Box{~}$}

\newenvironment{eproof}
               {\noindent {\textbf{Proof.}}}
               {\proofend\\}

\newenvironment{fproof}
               {\noindent {\emph{\textbf{Preuve.}}}}
               {\proofend\\}

\newenvironment{enumeratei}{
 \renewcommand{\theenumi}{\roman{enumi}}
 \renewcommand{\labelenumi}{(\theenumi)} 
 \begin{enumerate}}{\end{enumerate}}

%%% Local Variables: 
%%% mode: latex
%%% TeX-master: t
%%% TeX-master: t
%%% TeX-master: t
%%% TeX-master: t
%%% TeX-master: t
%%% TeX-master: t
%%% TeX-master: t
%%% TeX-master: t
%%% TeX-master: t
%%% End: 

\title{Cramer's estimate for the exponential functional of a Levy process}
\author{Olivier \sc{Mejane}
        \\ \emph{Laboratoire de Statistiques et Probabilités,}
        \\ \emph{Université Paul Sabatier,
118 route de Narbonne}
     \\ \emph{31062 Toulouse France}   
     \\ email: Olivier.Mejane@lsp.ups-tlse.fr}

\maketitle

\section{Introduction}
The aim of this paper is to study the asymptotic behavior of the
 exponential functional $A_{\infty}=\int_0^{\infty} e^{\xi_s} \,ds$,
where $(\xi_t)_{t \geq 0}$ is
 a Lévy process such that Cramér's condition holds, that is to say there
 exists $\chi >0$ such that ${\dE}(e^{\chi \xi_1})=1$.
The precise result will be, under others conditions on $(\xi_t)_{t \geq 0}$,
that the tail of $A_{\infty}$ is decreasing like $t^{-\chi}$ when 
$t \to \infty$.
\\

One method to understand  this result is to start from the analogous problem
in discrete  
time:
the random difference equation $Y_n=M_n \,Y_{n-1}+ 1$ where
$(M_n)_{n \in {\dN}}$ is a sequence of i.i.d real variables,
and whose solution $Y_n$, under certain additional hypothesis (see
\cite{vervaat} and \cite{kesten}), converges in distribution to  
$R:=\sum_{k=0}^{\infty} M_1 \ldots M_{k-1} $.
In our case, the Lévy process $(\xi_t)_{t \geq 0}$ will play the role
of the random walk $(S_n:=\sum_{k=1}^{n}\log |M_k|)_{ n \geq 1}$ and
$A_{\infty}$ the role of the 
limit variable $R$.
\\

Let us explain a little the analogy. In the discrete case, for
all stopping-time $N$ which is finite almost-surely we have the following
identity in law (see \cite{vervaat}, 
 lemma 1.2)~:
\begin{equation} \label{discrete-time}
R \stackrel{d}{=} M_1 \ldots M_N  \, R + R_N \, .
\end{equation}
where $(R_n)_{n\in \dN}$ stands for the sequence of partial sums 
\mbox{$(\sum_{k=0}^{n} M_1 \ldots M_{k-1})_{n\in \dN}$}.
In continuous-time this identity is still valid if we replace 
$(R_n)_{n\in \dN}$ by the process \mbox{$(A_t:=\int_0^t e^{\xi_s} \,ds)_{t\geq
    0}$}; 
indeed the lemma 6.2 in \cite{carmona} implies that for all stopping-time $T$
 which is finite almost-surely we have :
\begin{equation} \label{continuous-time}
A_{\infty} \stackrel{d}{=}e^{\xi_T} A_{\infty} +A_T \, .  
\end{equation}

Kesten, in the above quoted article, found the asymptotic behaviour of the
distribution of $R$ ~:  
$${\dP}(R >t) \sim C \, t^{-\chi} \quad \mbox{for some constant} \, C > 0 \,
,$$ 
and he noticed that the
 proof (in the one-dimensional-case, which interests us here) relies
 essentially on Cramér's estimate for the random walk $S_n$ : if $\chi > 0$
 satisfies the Cramér's condition  
${\dE}(M_1^{\chi})=1$, then~: 
\begin{equation} \label{Cramer_for_the_max}
{\dP}(\max(S_0,S_1,\ldots) >t) \sim C \, e^{-\chi t},
\end{equation}
 for some constant $C$ when $t \to +\infty$.
\\
But Bertoin and Doney ( see \cite{cramer} ) have proved that Cramér's estimate
extends to Lévy processes.
Their proof relies on the introduction of what they call the associated
 Lévy process $X^*$, whose exponent is $\Phi^*(\la)=\Phi(\la+\chi)$.
It amounts to make the change of probability defined by the martingale
$(e^{\chi \xi_t})_{t \ge 0}$, as we do in our proof of Proposition 
\ref{renouv}. 
They use a Wald identity for this process, namely~:
$$\dE(H_1^*)= \dE(X_1^*) \dE(\tau_1^*),$$
where $\tau^*$ is the inverse process of the local time $L^*$ of the
reflected process at the supremum
\\
$(S_t^*-\xi_t^*= sup_{s \le t} \xi_s^* -\xi_t^*, t \ge 0)$ and
$H$ the ascending ladder process: $H_t^*=S_{\tau_t}^*=\xi_{\tau_t}^*$.
For sake of completeness, we give a proof of this Wald identity in the
Annex. 
\\

In fact, instead of the equivalence mentioned in (\ref{Cramer_for_the_max}),
 we can easily obtain an upper bound, and
this suffices to show that $A_{\infty}$ has moments of all orders
 $\alpha < \chi$ where $\chi$ is again the non-negative root of ${\dE}(e^{\chi
   \xi_1})=1$ . 
We will give a proof of this result (in the third section), which has an
interest by itself and will be used to obtain the more precise result
concerning the tail of $A_{\infty}$, that we give now~:

\begin{ethm} \label{nouveau}
Let $\xi$ be a Lévy process, with Lévy exponent $\Phi$ (i.e
 ${\dE}(e^{\lambda \xi_t})=e^{-t \Phi(\lambda)}$)and fulfilling the following
Cramér's condition ~:

\begin{equation} \label {Cramer}  
\exists \chi >0 \mbox{ such that } \Phi(\chi)=0
\end{equation}

Notice that this can only happen if $-\infty \leq \mu := {\dE}(\xi_1) < 0$ .

We make besides the stronger hypothesis that
 $\Phi > -\infty $ on an interval $[0,\chi + \epsilon]$,
with $\epsilon > 0 $.

At last we assume that the law of $\xi_1$ is not arithmetic. 
\\
Then the exponential functional 
$A_{\infty} := \int_0^{\infty} e^{\xi_s} \,ds$ is well defined and there
exists some constant $C>0$ such that when $t \to +\infty$ : 

\begin{equation} \label {tail}  
{\dP} (A_{\infty}>  t) \sim C \, t^{-\chi} 
\end{equation}
\end{ethm}

The proof we give in the third part of the article is greatly inspired by
Goldie (\cite{goldie}) who found a simpler proof of Kesten's result, 
that extends to the continuous case as we will show here. 

\section{Examples}
\subsection{Brownian motion with drift}
If $\xi_t=\sigma  B_t + \nu \,t$ is a brownian motion with negative drift
 ( $\nu < 0$ ), then $\Phi(\la)=\la( \frac{\sigma^2}{2} \la +\nu)$, so it
 satisfies the Cramér's condition with:  
$$\chi=\frac{-2\nu}{\sigma^2} \, .$$
In fact, in this case, we know explicitly the law of the exponential functional
(see for instance \cite{carmona})~:

$$\int_0^{\infty} e^{\sigma B_s + \nu s} \,ds \stackrel{d}{=}
\frac{2}{\sigma^2 \gamma_{-2\nu/\sigma^2}} \, ,$$ 
where $\gamma_m$ denotes a gamma variable with index $m$.
This implies easily the asymptotic behaviour of $A_{\infty}$ given
by Theorem \ref{nouveau}.

\subsection{Compound Poisson process with drift}
Let us take $\xi_t= -t -\eta_t$ where $\eta$ is a compound Poisson process 
with Lévy measure $$\mu(dx)=(a+b-1)be^{bx} \,dx \, ,x<0 \, ,$$
with $0<a<1<a+b$.
Then for $\la<b$ we have $\Phi(\la)=\frac{\la}{b-\la} (1-a-\la)$,
therefore it satisfies the Cramér's condition with 
$$\chi =1-a \, .$$

Here we also know the law of the exponential functional:
$$\int_0^{\infty} e^{\xi_s} \,ds \stackrel{d}{=} \frac{1}{\beta_{1-a,a+b-1}}
\, ,$$ 
where $\beta_{a,b}$ is a beta variable with parameters $a$ and $b$, so once 
again we could find easily (\ref{tail}).
 
\subsection{Opposite of a stable subordinator with drift}
Here we consider $\xi_t=-S_t+at$ where $(S_t)_{t\geq 0}$ denotes a standard
stable subordinator with index $ 0< \alpha<1$, and $a$ is a positive real.
Then $\Phi(\la)=\la^{\alpha} - a\la$ for $\la \geq 0$ so that $\Phi(\chi)=0$
for $\chi=a^{1/(\alpha -1)}$.
Thus: 
$${\dP} (\int_0^{\infty} e^{-S_u+au}\, du >  t) \sim 
\frac{c}{t^{a^{1/(\alpha -1)}}} \, .$$
This example is of greater interest since here we cannot compute the law  
of the exponential functional.
\section{Moments of the exponential functional}

\begin{eprop} \label{moments}
If Cramér's condition (\ref{Cramer}) is satisfied, then
 $${\dE} \left( A_{\infty}^{\alpha} \right) < \infty \, , \quad  \forall 
 \, 0 \leq \alpha < \chi$$
\end{eprop}
\begin{eproof}
The proof relies on the following lemma :
\begin{elem}
Under the precedent hypothesis, if we note 
$S_{\infty}=\sup_{t \geq 0} \xi_t$, we have
$$\dE(e^{\alpha S_{\infty}}) < \infty \quad \forall 
 \, 0 \leq \alpha < \chi$$ 
\end{elem}
Let us write $S_t=\sup_{0 \leq s \leq t} \xi_s$ for each $ t \geq 0$.
Since $\dE(e^{\chi \xi_t})=1$, the process $(e^{\chi \xi_t})_{t \geq 0}$ is a 
nonnegative martingale, to which we can apply Doob's Submartingale
Inequality for fixed $t >0$ and $x \in \dR^{+}$ : 
$$x \dP(\sup_{0 \leq s \leq t} e^{\chi \xi_s} \geq x) \leq \dE(e^{\chi
  \xi_t})=1$$ 
We obtain that for fixed $a\in \dR^{+}$ and for all $t > 0$,  
$\dP(S_t \geq a) \leq e^{-\chi a}$.
Since $S_{\infty}$ is finite almost-surely thanks to the fact that
 $\xi_t \to -\infty$ when $t \to +\infty$, it follows that
$\dP(S_{\infty} > a) \leq e^{-\chi a}$, which concludes the proof of the
lemma  
$\quad \Box$

Now let us fix $0 \leq \alpha < \chi$ and introduce the Lévy process
$(\xi'_t=\xi_t+k t,t \geq 0)$ with $k > 0$
small enough to ensure that $\xi'$ has the same properties as 
$\xi$: precisely we assume that $k+\dE (\xi_1) < 0$. On the one
hand we then have that $\xi'_t \to -\infty$ when $t \to +\infty$ and on the
other hand the Lévy exponent $\Psi$ of $\xi'$ still has a unique zero  
$\chi' > 0$ , such that $0 \leq \alpha < \chi'$ is $k$ is taken small enough. 
\\
Now we notice that :
$$A_{\infty}=\int_0^{\infty} e^{\xi'_s -ks} \,ds 
\leq e^{S_{\infty}'} \int_0^{\infty}e^{-ks} \,ds  =\frac{1}{k}
e^{S_{\infty}'}$$ 
writing as above $S_{\infty}'=\sup_{t \geq 0} \xi'_t$.
We deduce that:
$$\dE(A_{\infty}^{\alpha}) \leq  \frac{1}{k^{\alpha}} \dE(e^{\alpha
  S_{\infty}'})$$ 
If suffices to apply the lemma to the Lévy process $(\xi'_t)_{t \geq 0}$
to end the proof.
\end{eproof}

\section{Proof of the theorem}

\subsection{First step}
To prove (\ref{tail}) , it suffices to prove that~:

\begin{equation} \label{but}
\int_{-\infty} ^ t e^{-(t-v)} r(v) \, dv \xrightarrow[ t \to +\infty]{}
 C \quad \mbox{where} \quad  r(v)=e^{\chi v} {\dP} (A_{\infty}>e^v)
\end{equation}
 
It's a consequence of lemma 9.3 of \cite{goldie}, that we quote here
 for sake of completeness: 

\begin{elem}
Let $k>0$ and $X$ be a real random variable. If 
$\int _0^t  u^{k}{\dP} (X>u) \, du \sim C \, t$ when $t \to \infty$~,
 then ${\dP} (X > t)  \sim C \, t^{-k}$ when $t \to +\infty$
\end{elem}

%Thus we have to prove that there exists a constant $0<C<+\infty$ such that:
%\\
%$e^{-t} \int_0^{e^t} u^{\chi}{\dP} (A_{\infty}>u) \, du \to C$  when
%$t \to +\infty$, or equivalently, by change of variable:

If we introduce the function $K(t)=e^{-t}$
 for $t>0$ and equal to $0$ for $t\leq 0$ ,
 we can write the left member of (\ref{but}) as the convolution between $r$
 and $K$, and we will denote it by $\tilde r$. 
More generally for all function $f$ , we will note
 $${\tilde f} (t) = f*K(t) =\int_{-\infty} ^ t e^{-(t-u)} f(u) \, du \, .$$

The key ingredient of the proof will be then a renewal theorem.

\subsection{Second step}
We are now going to write ${\tilde r}$ in the following form~:
\begin{equation} \label{convol}
{\tilde r} (t)=\tilde g *\nu_{n-1}(t) + {\tilde \delta}_n(t) \quad \forall
n \geq 1 \,
\end{equation}
with appropriate functions $g$ and ${\tilde \delta}_n$ and measure $\nu_{n-1}$.
\\
For this we need a few notations:
Let $(T_i)_{i\geq 1}$ be a sequence of i.i.d variables with exponential law of
parameter $1$, and independent of the Lévy process $(\xi_t)_{t \geq 0}$. 
 Let $\Theta _n=\sum_{i=1}^{n} T_i$ for $n \geq 1$.
 The process $(S_n:=\xi_{ \Theta _n})_{n \geq 0}$ , with $S_0=0$, is a random
walk. 
Let us note that~:
\begin{equation} \label{mean}
{\dE} (S_1)=\int_0^{\infty} {\dE} (\xi_t) e^{-t} \, dt
={\dE} (\xi_1)
\int_0^{\infty} t e^{-t} \, dt={\dE} (\xi_1) \in [-\infty,0) 
\end{equation} 
so that 
$S_n \xrightarrow[ t \to +\infty]{a.s} -\infty$.

\begin{elem}
For all $n \geq 1$ we have $r(t)=g*\nu_{n-1}(t) + \delta_n(t)$, with~:
\begin{itemize}
\item
$g(t)=e^{\chi t} \left( {\dP} (A>e^t)-{\dP} (M A>e^t) \right)$ ,
with $M$ (respectively $A$) a random variable distributed as $e^{S_1}$ 
(respectively $A_{\infty}$), independent of $(\xi_t)_{t \geq 0}$ and
 $(T_i)_{i\geq 1}$, 
and $M$ independent of $A$.
\item
$\nu_n(dt)=e^{\chi t} \sum_{k=0}^{n} {\dP}(S_k \in dt)$
\item
$\delta_n(t)=e^{\chi t} {\dP} (e^{S_n}A>e^t)$
with $A$ as above.
\end{itemize}
\end{elem}
\begin{eproof}
First, by the identity in law between $A$ and $A_{\infty}$,
${\dP} (A_{\infty}>e^t) =  {\dP} (A>e^t)$; then by a different
 way of writing, we obtain that for all fixed $n \geq 1$  ~: 
$$
\begin{array}{cc}
{\dP} (A>e^t) 
 & = 
\sum_{k=1}^{n} {\dP} (e^{S_{k-1}}A>e^t)-
{\dP} (e^{S_k}A>e^t) \quad  +  \, {\dP} (e^{S_n} A>e^t) \\ 
& = \sum_{k=1}^{n} {\dP} (e^{S_{k-1}}A>e^t)-
{\dP} (e^{S_{k-1}} M A >e^t) \quad + \, {\dP} (e^{S_n} A>e^t) ,
\end{array}
$$

the last equality resulting from the independence between $(M,A)$ and
$(S_n)_{n \geq 0}$, and from the fact that this process is a random walk. 
 
Thus we have ~: 
$$
\begin{array}{cc}
r(t) &=  \sum_{k=0}^{n-1} 
\int_{\dR} e^{\chi (t-u)}\left( 
{\dP} (A>e^{t-u})-
{\dP} ( M A >e^{t-u})
\right) e^{\chi u} {\dP}(S_k \in \, du)
+ e^{\chi t}{\dP} (e^{S_n} A>e^t) \\ 
& =  \int_{\dR} g(t-u) \nu_{n-1}(du) +\delta_n(t)
\end{array}
$$

and this ends the proof $\quad \Box$
\\

Since ${\tilde r}=r*K$, ${\tilde g}=g*K$ and
${\tilde \delta}_n = \delta_n*K$, 
this lemma  obviously implies (\ref{convol}).
\end{eproof}

\subsection{Third step}

We are now going to show that~:
\begin{equation} \label{lim1}
\forall t \quad {\tilde \delta}_n(t) \to 0  \quad \mbox {when} \quad n\to
+\infty 
\end{equation}
and then, that there is a renewal measure $\nu$ such that~:
\begin{equation} \label{lim2}
\forall t \quad \tilde g *\nu_n(t) \to  \tilde g *\nu(t) \quad \mbox {when}
\quad n\to +\infty 
\end{equation}

\begin{itemize}
\item
The first point is the easier one~:
\\
$e^{S_n} \xrightarrow[ n \to +\infty]{a.s} 0$, so for fixed $t$,
$\delta_n(t) \xrightarrow[ n \to +\infty]{} 0 $.
We conclude by dominated convergence~:
${\tilde \delta}_n(t)=\int_{-\infty}^t e^{-(t-u)} \delta_n(u) \, du$
but $ 0 \leq \delta_n(u) \leq e^{\chi u}$ and 
$ \int_{-\infty} ^ t e^{(\chi+1)u} du < \infty$ since $\chi >0$.
\\

\item
To establish (\ref{lim2}), we shall need the following proposition~:
\\
\begin{eprop} \label{renouv}
Let be $\nu(dt)=e^{\chi t} \sum_{k=0}^{\infty} {\dP}(S_k \in \, dt)$.
Then $\nu$ is the renewal measure associated to some random walk
 $(Y_i)_{i \geq 0}$  such that $0< m:={\dE} (Y_1)<+\infty$
and that the law  of $Y_1$ is not arithmetic.
\end{eprop}

\begin{eproof}
The proof relies on the following change of probability~:
if for all $t \geq 0 $ we denote by ${\mathcal F}_t$ the natural filtration
 of the process $(\xi_t)_{t \geq 0}$, since $(e^{\chi \xi_t})_{t \geq 0}$ is a 
strictly positive $({\mathcal F}_t)_{t \geq 0}$-martingale, 
we can define a probability ${\dQ}$ on $\, \vee_{t\geq 0} {\mathcal F}_t$
by the change of probability ${\dQ}=e^{\chi \xi_t} {\dP}$
on ${\mathcal F}_t$, in other words~:
 
\begin{equation} \label{chgt_proba}
{\dE}_{\dQ}(X)={\dE}_{\dP}(e^{\chi \xi_t} X) \quad 
\mbox{for all bounded and } {\mathcal F}_t \mbox{-measurable} \,X \, 
\end{equation}

We easily check that under ${\dQ}$ , $(\xi_t)_{t \geq 0}$ is still
a Lévy process, with Lévy exponent 
$\Phi_{\dQ}(\lambda)=\Phi(\lambda+ \chi)$.
Since $\Phi$ is concave, and $\Phi (0)=\Phi (\chi)=0$ 
and last that $\Phi '(0)=-\mu  > 0$,
 we have
$${\dE}_\dQ(\xi_1)=-\Phi_\dQ '(0)=-\Phi '(\chi) > 0 \, .$$
Let then $(Y_i)_{i \geq 0}$ be a sequence of i.i.d variables, whose common
distribution is the law of  $S_1$ under $\dQ$. 

Let us consider the renewal measure $U$ associated to this process 
$(Y_i)_{i \geq 0}$, precisely~:
$$U(dt)=\sum_{k=0}^{\infty} \dQ (S_k \in \, dt) \,$$

Since  $S_k=\xi_{ \Theta _k}$ is ${\mathcal F}_{\Theta_k}$-measurable,
 and since the formula (\ref{chgt_proba}) is still true with any almost-surely
 finite stopping-time  
$T$ instead of $t$, in particular for $T=\Theta_k$, we obtain~: 

 $${\dQ}(S_k \in dt)={\dE}_{\dP}(1_{\{\xi_{\Theta _k} \in dt\}} e^{\chi
   \xi_{\Theta _k}})=e^{\chi t}\dP(\xi_{\Theta _k} \in dt)$$ 
and thus,
$$U(dt)=e^{\chi t}\sum_{k=0}^{\infty} \dP(S_k \in \, dt) \, .$$
 As a consequence $U$ is exactly the measure $\nu$ introduced in the
 proposition.  
\\

What is left to be proved is that $(Y_i)_{i \geq 0}$ fulfills the hypothesis
of the proposition:
\\

-The sign of the mean is a result of the same calculus that in (\ref{mean}) ~: 
$m:= \dE (Y_1)=\dE_\dQ(\xi_1) > 0$. 
Moreover if we had $m=\infty$ that would be in contradiction with the
fact that $\Phi > -\infty$ on a neighborhood of $\chi$.
\\

-For the second point we first notice that if the law of $S_1$ is  
arithmetic under $\dQ$ then it is also the case under $\dP$ :
indeed if there is some $\lambda$ such that $\dQ(S_1 \in \lambda{\dZ})=1$,
that means that $\dE_\dP(1_{\{S_1 \in \lambda{\dZ}\} } e^{\chi S_1})=1$, but
the variable $e^{\chi S_1}$ being nonnegative and having a mean equal to $1$
under $\dP$ (since $\dE (e^{\chi \xi_{\Theta_1}} )=\int_0^{\infty} \dE
(e^{\chi \xi_t}) e^{-t} \, dt 
=\int_0^{\infty} e^{-t} \, dt= 1$),
this implies that $\dP(S_1 \in \lambda{\dZ})=1$.
But the law of $S_1=\xi_{\theta_1}$ cannot be arithmetic unless that of
$\xi_t$ is arithmetic for all $t$, which is excluded by hypothesis. 
\end{eproof}
\end{itemize}

To prove that ${\tilde g} *\nu_n(t) \to {\tilde g} *\nu(t)$, and then apply
a renewal theorem to $\nu$, we have now to show 
the direct Riemann-integrability of 
${\tilde g}$. This will be the fourth step, but for the moment
let us assume this result.
We then know that $|\tilde g|*\nu (t) < \infty$ for all $t$.
This means that~: 
$\dE(\sum_{k=0}^{\infty} e^{\chi S_k} | {\tilde g}(t-S_k)| ) < \infty$
and we deduce that 
$$\tilde g *\nu_n(t)=\sum_{k=0}^{n}\dE ( e^{\chi S_k} {\tilde g}(t-S_k))
\xrightarrow [ n \to + \infty]{} {\tilde g} *\nu(t)$$ 
 
We have thus proved the two points (\ref{lim1}) and (\ref{lim2}).
% ,and
% that, together with (\ref{etape1}) , implies ${\tilde r} (t) = {\tilde g}
 % %*\nu(t)$. 

%Now ${\tilde g}$ being directly Riemann-integrable, the proposition
% (\ref{renouv}) enables to apply the renewal theorem for the measure $\nu$~:
%$${\tilde g} *\nu(t)=\int_{\Bbb R} {\tilde g}(t-u) \nu(du) 
%\xrightarrow[t \to + \infty]{} 
%\frac{1}{m} \int_{\Bbb R} {\tilde g} \.$$
 
%So the proof is finished by taking 
%$C=\frac{1}{m} \int_{\Bbb R} {\tilde g}=\frac{1}{m} \int_{\Bbb R} g$,
% as soon as we have proved this last result~:

\subsection{Fourth step}
\begin{eprop}
${\tilde g}$ is directly Riemann-integrable.
\end{eprop}
\begin{eproof}
The key is the following lemma, the demonstration of which is given in
\cite{goldie} (p.143) ~: 

\begin{elem}
If $f \in L^1(\dR)$, then $\tilde f$ is directly Riemann-integrable.
\end{elem}

We are thus going to prove that   
$g(t)=e^{\chi t} \left( \dP (A>e^t)-\dP (M A>e^t) \right)$ is 
in $L^1(\dR)$.

Firstly, by the change of variables $u=e^t$,
$$\int_\dR |g(t)| \, dt =
\int_0^{\infty} u^{\chi-1} ( |\dP (A>u)-\dP (M A>u) )| \, du \,. $$

Now for any almost-surely finite stopping time $T$, by using the strong Markov 
property at time $T$ for the Lévy process $\xi$, we have the following
 decomposition of the exponential functional $A_{\infty}$~:
\begin{equation} \label{Markov}
A_{\infty}=A_T+e^{\xi_T} \int_0^{\infty} e^{\xi_{s+T} - \xi_T} \, ds
= A_T+e^{\xi_T}{\tilde A}_{\infty}^T \, ,
\end{equation}
the variable ${\tilde A}_{\infty}^T :=\int_0^{\infty} e^{\xi_{s+T}-\xi_T} \,
ds$ 
having the same law as $A_{\infty}$ and being independent of 
$(\xi_s)_{0 \leq s \leq T}$.
The identity (\ref{continuous-time}) mentioned in the introduction was a direct
consequence of (\ref{Markov}).
By using (\ref{Markov}) with $T=\Theta_1$ , we first obtain 
 $(M,A) \stackrel{d}{=} ( e^{\xi_{\Theta_1}}, {\tilde
   A}_{\infty}^{\Theta_1})$, and so $\dP (M A>u)=\dP
 (e^{\xi_{\Theta_1}}{\tilde A}_{\infty}^{\Theta_1}>u)$; then we write $\dP
 (A>u)=\dP (A_{\infty}>u)$ ,which eventually leads us to~: 
$$|\dP (A>u)-\dP (M A>u)|= \dP(e^{\xi_{\Theta_1}}{\tilde
  A}_{\infty}^{\Theta_1} \leq u < A_{\infty})$$ 

Hence~:
$$\int_\dR |g(t)| \, dt =\int_0^{\infty} u^{\chi-1} | \dP (A>u)-\dP (M A>u) |
\, du 
=\frac{1}{\chi} \dE \left( A_{\infty}^{\chi} -(e^{\xi_{\Theta_1}}{\tilde
    A}_{\infty}^{\Theta_1})^{\chi} \right)$$ 

Two cases have to be distinguished~:

\begin{itemize}
\item
{\it First case:}  $0<\chi\leq 1$
\\
Then the following inequality holds : $|x^{\chi}-y^{\chi}| \leq |x-y|^{\chi}$
for all nonnegative $x,y$ .
Thus $$\frac{1}{\chi} \dE \left( A_{\infty}^{\chi} -(e^{\xi_{\Theta_1}}{\tilde
    A}_{\infty}^{\Theta_1})^{\chi} \right) 
\leq \frac{1}{\chi} \dE(A_{\infty}-e^{\xi_{\Theta_1}}{\tilde
  A}_{\infty}^{\Theta_1})^{\chi}=\frac{1}{\chi} \dE (A_{\Theta_1} ^{\chi}) \,
.$$ 
To see that the last right-member is finite, we can use the same method that
in the proof of Proposition \ref{moments}; indeed $A_{\Theta_1}$ can be
 seen as the terminal value of the process $(\int_o^t e^{\xi_s^{(1)}}
 \,ds)_{t\geq 0}$ where $\xi^{(1)}$ is the initial Lévy process killed at the 
independent exponential time $\Theta_1$. His Lévy exponent is 
$\Phi^{(1)}=\Phi + 1$. Since $\Phi(\chi)=0$ and since there exist 
$\epsilon >0$ for which $\Phi > -\infty $ on $[0,\chi + \epsilon]$, by
continuity we can find $\la_0 > \chi $ such that  
$\Phi^{(1)}(\la_0) >0$.
Then if we now use the martingale~: 
$$(e^{\la_0 \xi_s^{(1)} + \Phi^{(1)}(\la_0) s)})_{t \geq 0} \, ,$$ we
deduce, as in Proposition \ref{moments}, that the
exponential functional associated to the killed Lévy process $\xi^{(1)}$
has moments of all order $\alpha < \la_0$; in particular with $\alpha=\chi$
we obtain exactly that $\dE (A_{\Theta_1} ^{\chi}) < \infty$.  
\\

\item
{\it Second case:}  $1<\chi$
\\
This time we use the inequality
$|x^{\chi}-y^{\chi}| \leq \chi |x-y| (\max(x,y))^{\chi-1}$, which leads to
 the upper bound~:

\begin{equation} \label{majo1}
\int_\dR |g(t)| \, dt \leq 
\dE \left( (A_{\infty}- e^{\xi_{\Theta_1}}{\tilde A}_{\infty}^{\Theta_1})
A_{\infty}^{\chi-1} \right) \, .
\end{equation}

The right-member can also be written
$$\dE\left( A_{\Theta_1} (A_{\Theta_1}+e^{\xi_{\Theta_1}}{\tilde
    A}_{\infty}^{\Theta_1})^{\chi-1} \right) \, . $$ 
Since $|x+y|^{\chi} \leq c_r ( |x|^{\chi}+|y|^{\chi})$ for some constant
$c_r$, we obtain~: 

\begin{equation} \label{majo2}
\int_\dR |g(t)| \, dt \leq 
 c_{\chi-1} \left[ \dE (A_{\Theta_1}^{\chi}) + 
\dE \left( A_{\Theta_1} ({\tilde A}_{\infty}^{\Theta_1}
e^{\xi_{\Theta_1}})^{\chi-1} \right) \right] \, .
\end{equation}

We have already seen that $\dE (A_{\Theta_1}^{\chi}) < \infty$.
Concerning the second term, we first use the independence between  ${\tilde
  A}_{\infty}^{\Theta_1}$ and $(\xi_s)_{0 \leq s \leq \Theta1}$ ~: 

$$\dE \left( A_{\Theta_1} ({\tilde A}_{\infty}^{\Theta_1}
  e^{\xi_{\Theta_1}})^{\chi-1} \right)= 
\dE \left( ({\tilde A}_{\infty}^{\Theta_1})^{\chi-1} \right)
\dE \left( A_{\Theta_1} e^{(\chi-1)\xi_{\Theta_1}}) \right)$$

But $1<\chi$ so by Hölder's inequality ~:
$$\dE \left( A_{\Theta_1} e^{(\chi-1)\xi_{\Theta_1}}) \right)
\leq \dE ( A_{\Theta_1}^{\chi})^{1/{\chi}} 
\dE (e^{\chi \xi_{\Theta_1}} )^{(\chi-1)/{\chi}} $$

We have already seen that 
$\dE (e^{\chi \xi_{\Theta_1}} )= 1$
so~:
$$ \dE \left( A_{\Theta_1} e^{(\chi-1)\xi_{\Theta_1}}) \right) < \infty$$

Lastly ~:
$$\dE \left( ({\tilde A}_{\infty}^{\Theta_1})^{\chi-1} \right)=
\dE \left( A_{\infty}^{\chi-1} \right) < \infty$$

because $A_{\infty}$ has finite moments of all order
$0\leq \alpha < \chi$ (see Proposition \ref{moments}).

In conclusion, in this second case we still have $\int_\dR |g(t)| \, dt <
+\infty \, .$
\end{itemize}
\end{eproof}

\subsection{Conclusion}
 The second and third step imply that~: 
$${\tilde r} (t) = {\tilde g} *\nu(t) \, .$$
Now the Proposition \ref{renouv} and the last step enable us to apply
the renewal theorem to $\nu$ :
$$\tilde{g} *\nu(t)=\int_\dR {\tilde g}(t-u) \nu(du) 
\xrightarrow[t \to + \infty]{} 
\frac{1}{m} \int_\dR {\tilde g} \, .$$
So the proof is finished by taking 
$C=\frac{1}{m} \int_\dR {\tilde g}=\frac{1}{m} \int_\dR g$.

\subsection{Another method}

There exists a shorter proof of Theorem \ref{nouveau}, that we are 
going to detail here.
\\
Nevertheless, the previous proof has an interest by itself since
it shows that techniques used in discrete time can be adapted 
to continuous time.  
\\
Let us now explain this other method.
\\
The starting point is to notice that $A_{\infty}$ satisfies the random
difference equation~: 
$$A_{\infty}\, =M A'_{\infty}+Q$$
where on the right hand side
$$A'_{\infty}=\int_{0}^{\infty}\exp(\xi_{1+s}-\xi_1)ds\quad ,\quad
M=e^{\xi_1}\quad \hbox{and} \quad Q=\int_{0}^{1}e^{\xi_s}ds\,.$$
$A'_{\infty}$ is distributed as $A_{\infty}$ and is
independent of the pair $(M,Q)$, which enables to recover the
relation (\ref{continuous-time}) of the introduction when $T=1$.
\\
The idea is then to show that $M$ and $Q$ satisfy the
conditions of Kesten's Theorem (cf. e.g. Theorem 4.1 in \cite{goldie}),
which gives the conclusion.
\begin{itemize}
\item
The first condition $\dE(\ABS{M}^{\chi})=1$ is of course satisfied.
\item
The second one is $\dE(\ABS{M}^{\chi} \log^+\ABS{M}) < \infty .$
\\
This is true here thanks to the hypothesis that $\Phi > -\infty $ on an
interval $[0,\chi + \ep]$ for some $\ep > 0 $.
Indeed one has:
\begin{eqnarray*}
\dE(\ABS{M}^{\chi} \log^+\ABS{M})&=&\dE(e^{\chi \xi_1} \xi_1^+) \\
&\le &\frac{1}{\ep}\dE \PAR{ e^{\chi \xi_1} e^{\ep \xi_1}}  \\
&=&
\frac{1}{\ep} e^{-\Phi(\chi+\ep)} < \infty.
\end{eqnarray*}

\item
The last thing to check is that $\dE \ABS{Q}^{\chi} < \infty \, .$
\\
In fact one can prove that $Q^{\chi} \in L^{1+\ep_0}$ with
$\ep_0=\frac{\ep}{\chi}$ 
for all $\ep$ such that $\Phi(\chi+\ep) > - \infty$.
Indeed one first observes that $$0 \le Q^{\chi} \le \sup_{0 \le s \le 1}
e^{\chi \xi_s} \, .$$
But $(e^{\chi \xi_s})_{0 \le s \le 1}$ is a martingale bounded in 
$L^{1+\ep_0}$, since for all $0 \le s \le 1$:
$$\dE \PAR{e^{\chi(1+\ep_0) \xi_s}} = e^{-s \Phi(\chi+\ep)} \le 
e^{-\Phi(\chi+\ep)} < \infty \, .$$
Thus by Doob's Inequality in $L^{1+\ep_0}$, one concludes that:
$$\dE \PAR{ Q^{\chi(1+\ep_0)}} < \infty \, .$$
\end{itemize}

\section{Annex: proof of Wald identity for Lévy processes}
We first need some notations before stating that we call Wald identity,
 since it generalizes the so-called result for random
walks.
\\
Let $X$ be a Lévy process started at $0$, with Lévy exponent $\Psi$, such that
$0 < \mu\:=\dE(X_1) < \infty$. Let us write $S_t=\sup_{s \le t} X_s$;
let $L$ be a local time at $0$ for the reflected process at the supremum
 $S-X$ and $\tau$ its inverse process:
$\tau_t = \inf \{s > 0, L_s >t \}$.
Since $\mu\:=\dE(X_1) >0$,  $L_{\infty }= + \infty$ a.s and $\tau_t < \infty$
for all $t \ge 0$. 
Then one can define the ascending ladder height process $H$ by 
$H_t=S_{\tau_t}=X_{\tau_t}$ for all $t \ge 0$.

\begin{eprop} \label{Wald-identity}
With the previous notations, if there exists $K<0$ such that
\\
$0 < \Psi(K) < \infty$, then~:
\begin{equation} \label{Wald-relation}
\dE(H_1)= \mu \, \dE(\tau_1) \, .
\end{equation}
\end{eprop}

{\it Remark:}
This result applies to the Lévy process $X^*$ defined in the introduction.
Indeed $\dE(X_1^*)=-\Phi'(\chi) \in \dR_+^*$ (cf proof of Proposition
\ref{renouv}) and $\Phi^*(\la)=\Phi(\la+\chi) \in  (0,+\infty)$
for all $-\chi < \la < 0$.
  
\begin{eproof}

We first prove that $\tau_1$, which is a stopping-time with respect to the
filtration 
 $\cF_t =\si\BRA{ X_s, 0 \le s \le t}$, is integrable: 

\begin{elem}
Under the only assumption that $0 < \dE(X_1) \leq \infty$, one has
$\dE(\tau_1) < \infty$ .
\end{elem}
\begin{eproof}
Let $K$ be the Lévy exponent of the subordinator $\tau$.
 One knows (see for example \cite{bouquinbertoin}, p.166) that :
$$K(\la)=\exp\PAR{\int_0^{\infty} (e^{-t} - e^{-\la t})t^{-1} \dP(X_t \geq 0)
  \,dt} 
\, .$$
By writing for all $\la >0$ that $\ln \la=\int_0^{\infty} (e^{-t} - e^{-\la
  t})t^{-1} \, dt$, one obtains the following expression~:
$$\frac{K(\la)}{\la}=\exp\PAR{\int_0^{\infty} (e^{-\la t} - e^{-t})t^{-1}
  \dP(X_t < 0) \,dt}\, .$$ 
But, since $\lim_{t\to +\infty} X_t = +\infty$,  
$\int_1^{\infty} t^{-1} \dP(X_t < 0) \,dt < \infty$ , and by dominated
 convergence one concludes that~:
$$K'(0)=\lim_{\la \to 0} \frac{K(\la)}{\la} = 
\exp\PAR{\int_0^{\infty} (1 - e^{-t})t^{-1} \dP(X_t < 0) \,dt}< \infty \quad
\Box $$ 
\end{eproof} 

Then we prove the integrability of the terminal value of the infimum process
\\
$(I_t=\inf_{s \le t} X_s, t \ge 0)$:

\begin{elem}
Let us write $I_{\infty}=\inf_{t \geq 0} X_t$.
Then $\dE (I_{\infty}) > -\infty \, .$ 
\end{elem}
\begin{eproof}
First we notice that $I_{\infty}$ is finite a.s. since $\mu > 0$.
\\
We introduce the dual process 
${\tilde X}=-X$, whose Lévy exponent is 
${\tilde \Psi}(\lambda)= \Psi(-\lambda)$.
\\
Writing $M_t=\sup_{s \leq t}
{\tilde X}_s$ for $t\geq 0$, we are going to prove that $M_{\infty}$ has
exponential moments of order $\alpha$ for all $0 \leq \alpha < -K$, which a
fortiori implies that $M_{\infty}$ is integrable, and the conclusion will
follow since $I_{\infty}=-M_{\infty}$.
\\
Since ${\tilde \Psi}(-K)=\Psi(K) < \infty$, the process $(e^{-K {\tilde X}_t +
  {\tilde \Psi}(-K) t)})_{t \geq 0}$ is well defined and is a
non-negative $(\cF_t)_{t \ge 0}$-martingale , so if we fix $t >0$ and $x \in
\dR^{+}$, we  
have by Doob's Submartingale Inequality~:
$$x \, \dP(\sup_{0 \leq s \leq t} e^{-K {\tilde X}_s +{\tilde \Psi}(-K)s} \geq
x) \leq \dE(e^{-K {\tilde X}_t +{\tilde \Psi}(-K) t})=1 $$

Since ${\tilde \Psi}(-K) >0$, this implies that $x \, \dP(\sup_{0 \leq s \leq
  t} e^{-K {\tilde X}_s} \geq x) \leq 1,$ 
 and thus that for all fixed $a \in \dR^{+}$ and all $t > 0$, $\dP(M_t \geq a)
 \leq e^{Ka}$. 
 It follows that 
$\dP(M_{\infty} > a) \le e^{Ka}$ 
which ends the proof of the lemma.
 $\quad \Box$
\end{eproof}

We are now able to prove (\ref{Wald-relation}).
Using that $(X_t - \mu \, t, t \geq 0)$ is a $(\cF_t)_{t \ge 0}$-martingale,
we deduce that 
 for all $n \ge 0$, 
\begin{equation} \label{stopping}
\dE(X_{\tau_1 \wedge n}) = \mu \, \dE(\tau_1\wedge n)
\end{equation} 
Writing $X_t=X_t^+ - X_t^-$ with 
$X_t^+ =\max(X_t,0)$ and $X_t^-=\max(-X_t,0)$, we have
\\
$X_{\tau_1 \wedge n}^+ \leq X_{\tau_1\wedge n} - I_{\infty}$,
 hence using (\ref{stopping}) and Fatou's lemma, one obtains~:
$$\dE(X_{\tau_1}^+) \leq \liminf_{n \to +\infty} \dE(\mu \, (\tau_1 \wedge n)
- I_{\infty}).$$ 
Thanks to the previous lemmas, 
$\liminf_{n \to +\infty} \dE(\mu \, (\tau_1 \wedge n) - I_{\infty})= \mu \,
\dE(\tau_1)-\dE(I_{\infty}) < \infty$, thus 
$$0 \leq \dE(H_1)=\dE(X_{\tau_1}) \leq \dE(X_{\tau_1}^+)< \infty \, .$$

Since $I_{\infty} \leq X_{\tau_1 \wedge n} \leq S_{\tau_1}=H_1$,
one concludes from (\ref{stopping})
 by dominated convergence that
$$ \dE(X_{\tau_1})= \mu \, \dE(\tau_1). $$

\end{eproof}

\bibliographystyle{plain}
\bibliography{biblio}

\end{document}